\documentclass[conference]{IEEEtran}

\usepackage{amssymb}

\usepackage{multirow}
\usepackage{cite}
\usepackage {siunitx}

\usepackage{float}

\ifCLASSINFOpdf
 \usepackage[pdftex]{graphicx}

\else
\fi
\usepackage[cmex10]{amsmath}
\begin{document}

\title{Exploring Wake Interaction for Frequency Control in Wind Farms}

%\author{
%    \IEEEauthorblockN{Hesamoddin Marzooghi\IEEEauthorrefmark{1}\IEEEauthorrefmark{$,1$}, \emph{Student MIEEE}, David J. Hill\IEEEauthorrefmark{2}\IEEEauthorrefmark{$,1,2$}, \emph{Fellow MIEEE} and Gregor Verbi\v{c}\IEEEauthorrefmark{3}\IEEEauthorrefmark{$,1$}, \emph{Senior MIEEE}\\}
%    \IEEEauthorblockA{\IEEEauthorrefmark{$1$}School of Electrical and Information Engineering, The University of Sydney, Sydney, New South Wales, Australia\\}
%    \IEEEauthorblockA{\IEEEauthorrefmark{$2$}Department of Electrical and Electronic Engineering, The University of Hong Kong, Hong Kong\\}
%    \IEEEauthorblockA{\IEEEauthorrefmark{1}hesamoddin.marzooghi@sydney.edu.au, \IEEEauthorrefmark{2}david.hill@sydney.edu.au, %\IEEEauthorrefmark{3}gregor.verbic@sydney.edu.au\\}
  
%}

\author{\IEEEauthorblockN{Ahmad Shabir Ahmadyar\IEEEauthorrefmark {1}  and Gregor Verbi\v{c}\IEEEauthorrefmark {1}}
\IEEEauthorblockA{\IEEEauthorrefmark {1}School of Electrical and Information Engineering, The University of Sydney, Sydney, Australia\\}
    \IEEEauthorblockA{Emails:\{ahmad.ahmadyar, gregor.verbic\}@sydney.edu.au}}

%Telephone: (800) 555--1212, Fax: (888) 555--1212}


\IEEEpeerreviewmaketitle

\maketitle

%\doublespacing

\begin{abstract}
The increasing integration of wind generation is accompanied with a growing concern about secure and reliable power system operation. Due to the intermittent nature of wind, the base-load units need to cycle significantly more than they were designed for, resulting in reduced life cycle and increased costs. Therefore, it is becoming necessary for wind turbines to take part in frequency control, and reduce the need for additional ancillary services provided by conventional generators. In this paper, we propose an optimised operation strategy for the wind farms. In this strategy, we maximise the kinetic energy of wind turbines by an optimal combination of the rotor speed and the pitch angle. We exploit the wake interaction in a wind farm, and de-load some of the up-wind turbines. We show that the kinetic energy accumulated in the rotating masses of the WTs can be increased compared to the base case without compromising efficiency of the wind farm. In a specific system, we show that by implementing this strategy, and injecting the stored kinetic energy of the WTs' rotors into the system during a frequency dip, we can delay the system frequency nadir up to \SI{30}{\second}.   
\end{abstract}


% \fi

\section{Introduction}
\PARstart The dwindling fossil fuel resources, and their associated greenhouse gas emissions that significantly contribute to global warming, are encouraging nations to move toward renewable energy sources (RES). Among the RES, wind is one of the most economically viable options. In 2012, with 19\% growth, wind's annual growth was the highest among all RES\cite{GlobalWindEnergyCouncil2013}. USA, the world's largest electricity consumer, is planning to produce 20\% of its electricity from wind by 2030 \cite{InternationalInstituteforAppliedSystemsAnalysis2012}. Denmark, one of the pioneers of wind technology and the country with the highest penetration of wind in its electricity system, set the target of achieving 50\% of its electricity from wind power by 2020, and 100\% renewable by 2035 \cite{GlobalWindEnergyCouncil2013}. In Australia, wind is anticipated to play a major role in reaching the 20\% renewable energy target (RET) by 2020\cite{GlobalWindEnergyCouncil2013}. It has been predicted that from 2014 to 2035, under the 450ppm scenario, US\$3,027 billion will be invested in wind generation technologies, which is the highest investment among all renewable and conventional generation technologies \cite{IEA2014}. It appears that penetration of wind energy in the power system of most countries will keep increasing in the foreseeable future\cite{GlobalWindEnergyCouncil2013}. 
\par Although financially wind is one of the most viable RES, and can compete even with some of the conventional resources \cite{IEA2012}, technically, it is considered a less reliable resource because of its intermittent nature. Nowadays, variable speed wind turbines (VSWT), including Type III, using doubly fed induction generators (DFIG), and Type IV, with a generator connected through a fully rated converter (FRC), are considered the most promising technologies because of their ability to optimise the power extraction for over a wide wind speed range, as well as for being able to comply with grid codes' connection requirements. Therefore, approximately 95\% of wind turbines (WTs) installed all around the world are either Type III or IV\cite{BABUN.RAMESH2013}. Nonetheless, compared to the fixed speed wind turbines (FSWT), VSWTs do not have inherent inertial response \cite{Mullane2005}. A study in \cite{Doherty2010} has suggested that high penetration of Type III WTs can change the system frequency behaviour, which is characterized by the rate of change of frequency (ROCOF) and frequency nadir, so it is necessary for system operators to address these issues. Therefore, to increase integration of wind generation in power systems, wind will need to offer frequency control ancillary services, which will in turn reduce the pressure on conventional generators \cite{Doherty2010,Miller2010,Troy2010}. According to \cite{Ping-KwanKeung2009}, kinetic energy released by a WT exceeds that released by a synchronous generator (SG), which makes wind generation even more attractive to use for inertial contribution, including primary frequency control. The wind speed namely doesn't change significantly in short durations while the primary control is active, so the kinetic energy stored in rotors of partially loaded WTs can be used in frequency regulation. Due to the effective decoupling between the mechanical and the electrical systems of a WT, this capability needs to be emulated through appropriate control. 

\par An obvious issue in participation of WTs in frequency control is that it requires a WT to operate below its optimum power output for a period of time, which negatively impacts the efficiency of the WT. To partly overcome the negative financial impacts associated with spilling the wind energy, we propose to take advantage of the wake interaction within a wind farm (WF). The extraction of energy from wind by a WT namely results in a disturbed wind flow behind the WT, which can cause fatigue to the down-WTs, and in turn increases the maintenance cost as well as shortens the WT's life-cycle. An interesting approach how to deal with these issues has been proposed in \cite{Madjidian2011}. Using a stationary wake model, it has been shown that partially operating up-WTs not only reduces turbulence levels, but also improves the row efficiency of the WF. As a result, partial de-loading of the up-WTs in a WF not only offers frequency control ancillary services, but also reduces turbulence levels for the down-WTs.
 
\par  Therefore, we propose an optimised operation strategy for a WF, in which we de-load some up-WTs to maximise the overall kinetic energy of the WF. We do this by optimising the rotor speed $\omega$ and the pitch angle $\beta$ of a WT. We implement the optimised control approach initially proposed in \cite{Zertek2012}, and the wake model developed in \cite{Madjidian2011}. For particular de-loadings, we show that not only can the kinetic energy accumulated in the rotating masses of the WTs be increased, but also the  overall output power of the WF does not change significantly in a wide range of wind speeds compared to the base case. We cast the problem as a constrained non-linear optimisation problem, and solve it using the pattern search algorithm. Although the optimisation problem has multiple optima, we show that a good-quality solution can be found which can be readily implemented in a control algorithm.
\par The paper is organised as follows. In Section II, we briefly review the conventional strategies for participation of wind power in frequency control. Wake models are introduced in Section III. Section IV presents the proposed optimised operation strategy for frequency control of a WF. In Section V, the proposed control approach is implemented in a simple test system and results are evaluated. Section VI concludes the paper.
	
\section{Participation of wind power in frequency control}
\par
\begin{figure}
\centering
\includegraphics [width=9cm, height=6.8cm]{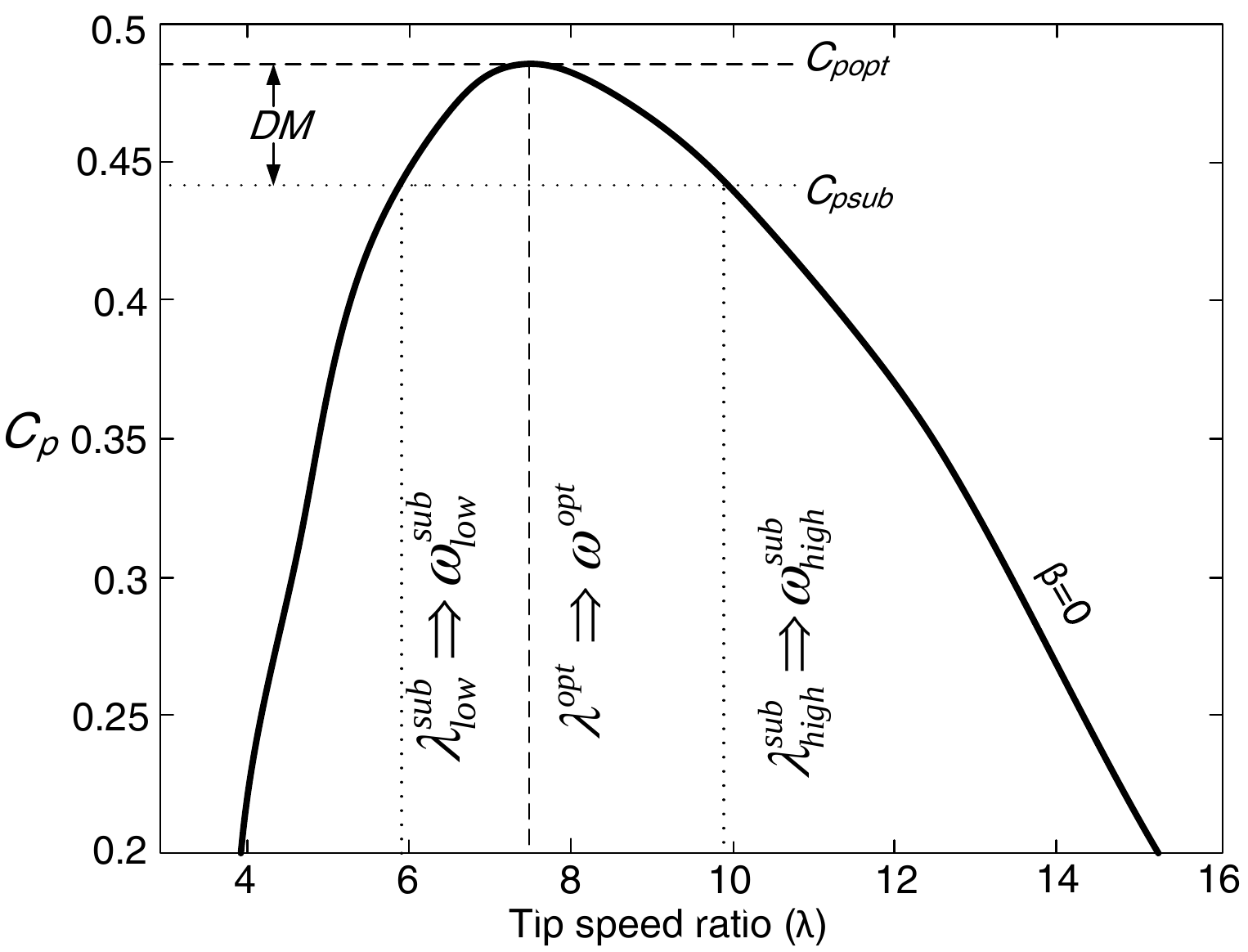} 
\caption{Power coefficient characteristic and de-Loading Strategy of Type III WT.}
\label{Figure_1}
\end{figure}

 There are mainly two options for control of frequency by WTs: \emph{inertial response} and \emph{de-loaded operation}\cite{Zertekb2012}. Since the stator and the rotor of a VSWT are decoupled by the power electronic converter, an additional control loop is required to make the WT inertia available to the system, which is referred to as synthetic inertia in some references\cite{Morren2006}. The addition of this new loop can provide up to 20\% extra power to the system for up to \SI{10}{\second} during a frequency dip \cite{Tarnowski2009}.  Nonetheless, this inertial response decreases the rotor speed which consequently reduces the coefficient of performance. To recover the coefficient of performance, the kinetic energy of the rotor should be restored, which can result in another frequency event\cite{Tarnowski2009}. Unlike synthetic inertia where no wind is spilled, in de-loading strategy a WT needs to be permanently de-loaded for frequency control, and operate with a lower coefficient of performance\cite{Holdsworth2004}. Between the cut-in wind speed and the wind speed where the rotor speed reaches its maximum value, de-loading can be achieved by changing the rotor speed in proportion to the de-loading margin ($DM$). Once the maximum rotor speed is reached, de-loading is possible only by changing the pitch angle. Fig.~\ref{Figure_1} shows the performance coefficient $C_{p}$ as a function of the tip speed ratio $\lambda$ defined as $\lambda=\frac{R\omega}{v}$, where $R$ is the radius of the WT blade, $\omega$ is the rotor speed of the WT and $v$ is wind speed. This strategy is valid only for those wind speeds where $\omega<\omega^{max}$. As shown in Fig.~\ref{Figure_1}, to reduce $C_{p}$, we have to either increase $\lambda$ or decrease it. The aim is to increase the kinetic energy; therefore, we have to select $\lambda^{sub}_{high}$ which corresponds to $\omega^{sub}_{high} > \omega^{opt}$. There is another de-loading method which is based on the pitch control only in all wind speeds $v$ where $ v_{cut-in} <v<v_{cut-out}$.  In this paper, the optimised control approach initially proposed in \cite{Zertek2012}, which uses both $\beta$ and $\omega$ to de-load a VSWT is implemented.              
\section {Wake models}
\par Accurate modelling of wake effects in a WF is a challenging task. Some of the factors which are considered for wake modelling are: the distance between the WTs, the radius of the WTs, the geography of the site and the operating points of the WTs. Wake models fall into two broad categories: (i) experimental and (ii) analytical. \emph{Experimental wake models} are based on measurements from the WFs, and are specific to these WFs. \emph{Analytical wake models} are based on the laws of fluid dynamics, and are mainly classified into three subclasses \cite{Crespo1999}: (i) kinematic models, (ii) field models and (iii) roughness element models. \emph{Kinematic models} are based on conservation of momentum and start by modelling a single wake for a WF. Although these models are suitable for large WFs, some of them assume a constant value for the coefficient of thrust $C_T$, which makes them unreliable for accurate wake models because $C_T$ of a WT changes in every operation points\cite{Madjidian2011}. Some well-known kinematic models are: Jensen's model \cite{Jensen1983} and Frandsen's model \cite{Frandsen2007}. \emph{Field models} give wind speed at every point behind a WT \cite{Crespo1999}, so computationally they are more complicated to implement. \emph{Roughness element models} are subdivided into infinite cluster models and finite cluster models. In infinite cluster models, the WF is considered as a single element, and the effect of individual WTs are lost. Whereas, the finite cluster models give the wind speed on each row of the WTs. Most of these wake models are not suitable for control purpose because either they are too complicated or unreliable. A \emph{Stationary wake model} which was recently developed for control purposes is thus used in this paper. This wake model requires minimum data, and its parameters have clear representation, which makes it suitable for control purposes\cite{Madjidian2011}.
\begin{figure} 
\centering
\includegraphics [width=8cm, height=2cm]{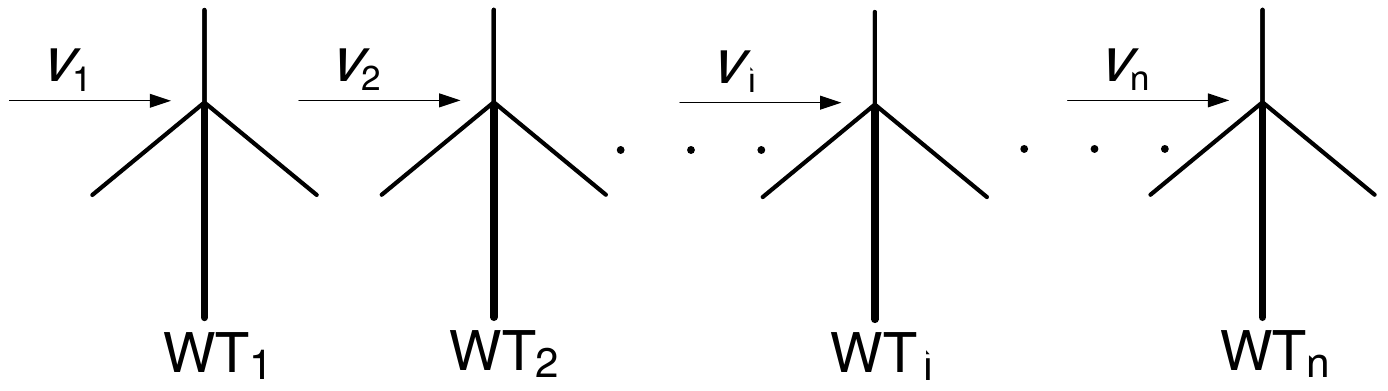} 
\caption{Row of WTs in a WF.}
\label{Figure_2}
\end{figure}
\par The stationary wake model is a suitable interaction model for a single row of WTs with the wind speed parallel to the row of the WTs. It maps the coefficient of thrust, which is directly linked to the wind speed deficit, the wind speed and the turbulence level of the up-WTs to identify their effect on the down-WTs\cite{Madjidian2011}. Considering the configuration of Fig.~\ref{Figure_2}, we can calculate the wind speed of $WT_{i+1}$:
\begin{align}
v_{i+1}=v_i + k^{'}(v_1 - v_i) - kv_1C_{Ti}
\end{align}  
where $0<k^{'}<1$ corresponds to the recovered wind speed, $0<k<k^{'}$ accounts the effect of the previous WT and $0<C_T(\lambda,\beta)<1$ is the thrust coefficient of the nearest up-WT. Distance parameters $k^{'}$ and $k$ are selected based on actual WF data. In the simplest approach, we need the data from three WTs. $k$ can be set based on $C_{T1}$ and $v_{2}$, and then $k^{'}$ can be set by having $k$, $C_{T2}$ and $v_{3}$. It is assumed that $k^{'}=0.35$ and $k=0.1$. The coefficient of performance $C_{p}(\lambda,\beta)$ and the coefficient of thrust $C_{T}(\lambda,\beta)$ are defined as:             
\begin{align}
C_p={\frac{2P_{mech}}{{\rho}{\pi}R^{2}v^{3}}}\\
C_T={\frac{2T_{F}}{{\rho}{\pi}R^{2}v^{2}}}
\end{align}
where $P_{mech}$ is the wind power transferred into mechanical power in the WT's rotor, $\rho$ is the air density, $R$ is  the radius of the WT's rotor, wind speed is $v$ and $T_{F}$ is the rotor thrust. Considering (2) and (3), and $P=\omega T$ a direct relationship between $C_{p}(\lambda,\beta)$ and $C_{T}(\lambda,\beta)$ in all operation regions where $C_{T}<{\frac{8}{9}}$ can be derived as follows \cite{Hasen2008}:  
\begin{align}
&C_p={\frac{1}{2}}(1 + \sqrt{1-C_T})C_{T}
\end{align}
where $C_p$ can either be given in a look-up table \cite{Jonkman2009}, or it can be approximated using curve-fitting\cite{Ackermann2005}. The above equations show that any changes in $C_{p}(\lambda,\beta)$ cause $C_{T}(\lambda,\beta)$ to change as well. Reducing $C_{p}(\lambda,\beta)$ for de-loaded operation results in a lower value of $C_{T}(\lambda,\beta)$, which results in less energy being extracted from the wind, which in turn increases the wind downstream and reduces the turbulence.

\section{Optimised operation strategy of wind farms for frequency control}

\subsection{Frequency control strategy}

\par
\begin{figure}
\centering
\includegraphics [width=8.7cm, height=6.5cm]{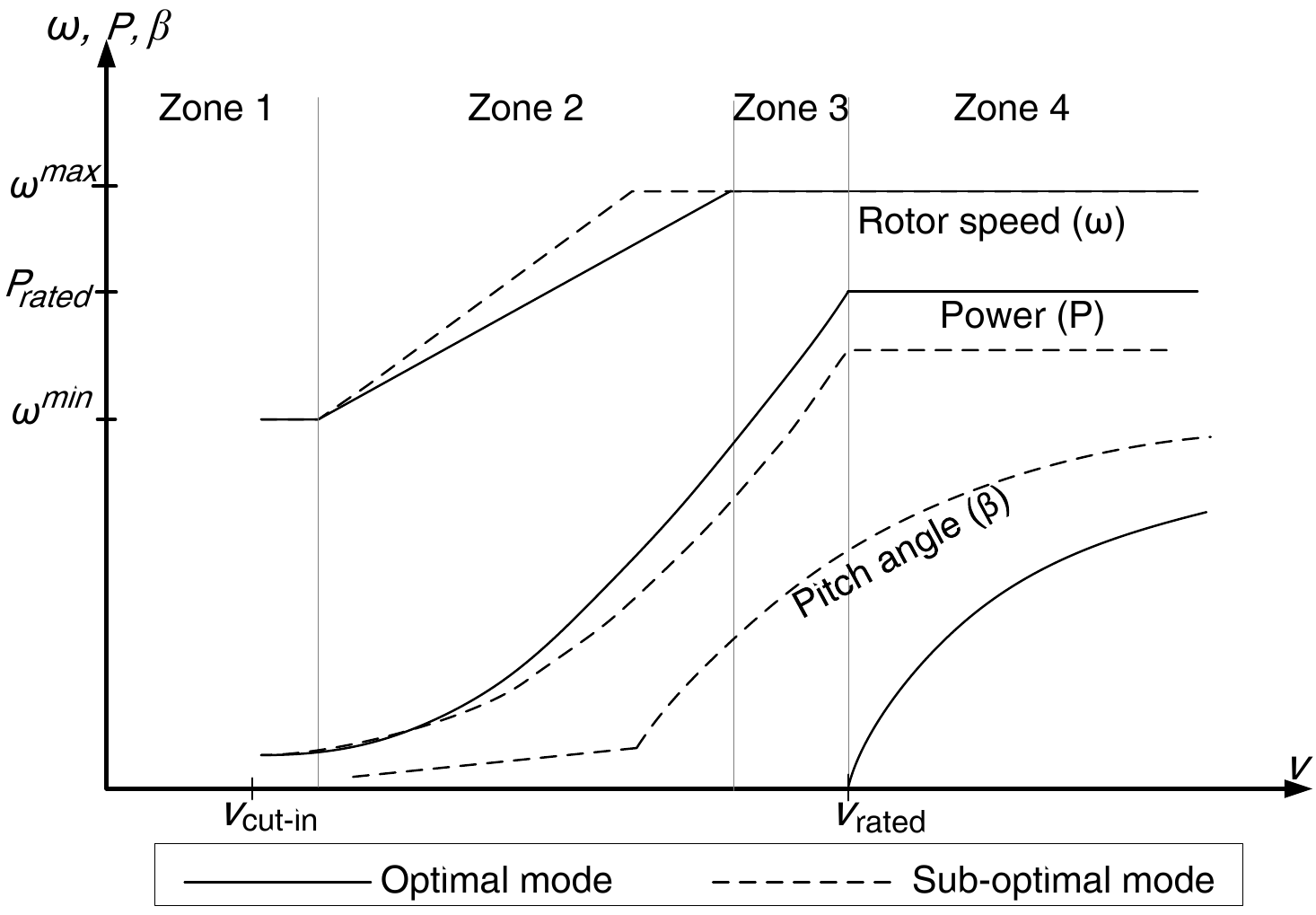} 
\caption{Power, rotor speed and pitch angle characteristic of a Type III WT.}
\label{Figure_3}
\end{figure}

 Fig.~\ref{Figure_3} illustrates the rotor speed, the pitch angle and the power characteristic of a Type III WT under the optimal and the sub-optimal (de-loaded) operation mode in four different zones. Under the optimal operation mode, the rotor speed is constant in Zones 1 and 3, and $C_p$ is not optimal. In Zone 2, any changes in the wind speed cause the rotor speed to change in order to maximise $C_p$. In Zone 4, the rotor speed is at its limit. Therefore, the pitch control operates to limit input the mechanical power to its rated value.      
\par Under the optimised sub-optimal operation mode\cite{Zertek2012} in Zones 1 and 2, both the rotor speed $\omega$ and the pitch angle $\beta$ vary to optimise the kinetic energy of the rotor, defined as:
\begin{align}
E_{k}=H{\omega}^{2}
\end{align}     
where $H$ is the normalized inertia of the WT, and $\omega$ is the rotor speed in the sub-optimal operation mode.

\par In Zones 3 and 4, since the rotor speed is maximum $\omega^{max}$, the only control variable is $\beta$, and no additional kinetic energy can be achieved. In Zone 2 under the sub-optimal mode, $\omega_{high}^{sub} > \omega^{opt}$, so Zone 2 becomes narrower as the $DM$ increases, and this limits the available $DM$\cite{Zertek2012}.

\begin{figure*}
\centering
\includegraphics [width=18cm, height=8.5cm]{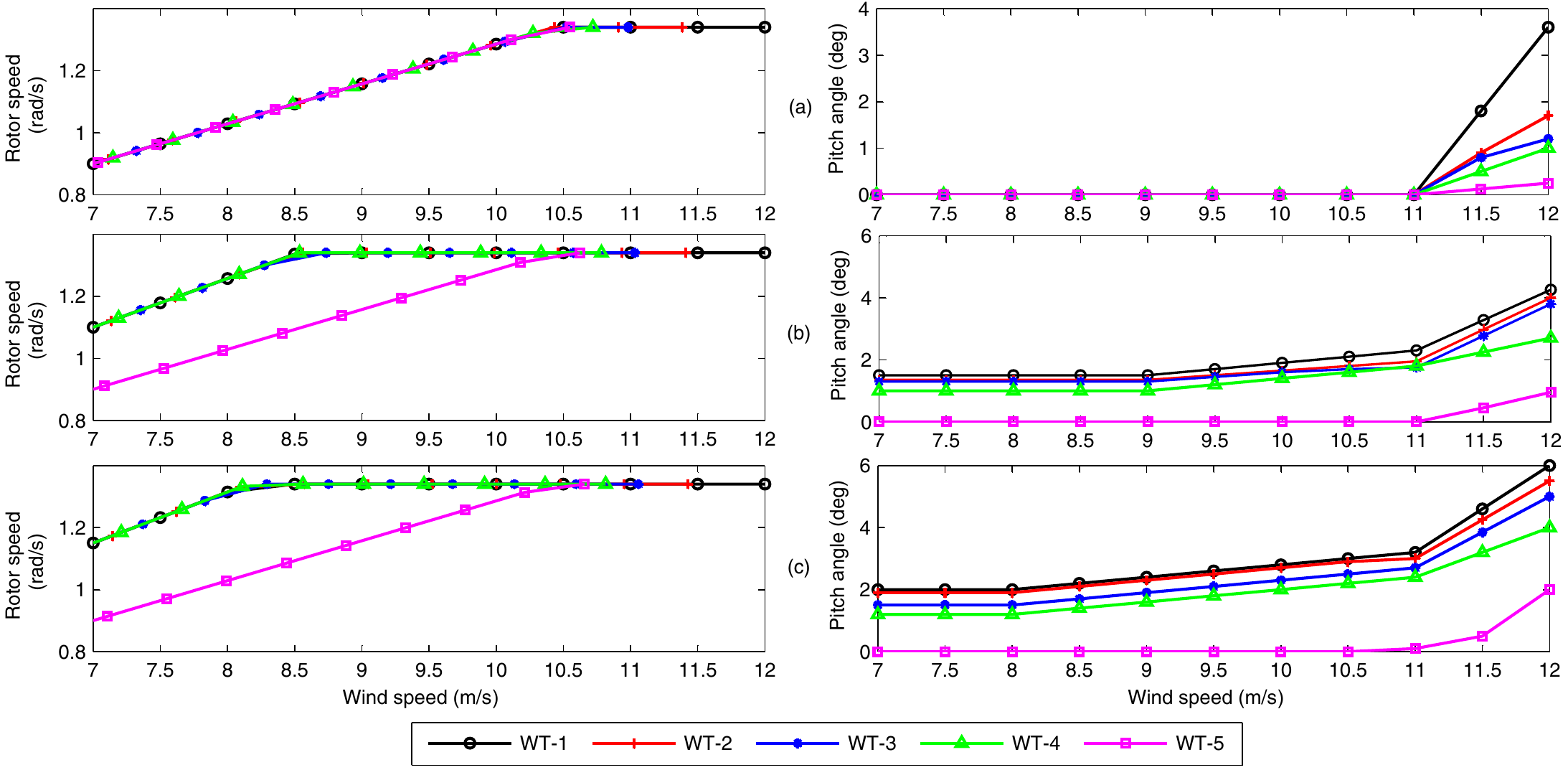} 
\caption{Pitch angle and rotor speed of the Type III WTs. a) $DM=0\%$, b) $DM=5\%$, c) $DM=10\%$.}
\label{Figure_4}
\end{figure*}
 
\subsection{The optimisation problem}

\par In this paper, we optimise the total kinetic energy of a WF by de-loading some of the up-WTs using a combination of the rotor speed $\omega$ and the pitch angle $\beta$. We consider a WF with a single row of identical Type III WTs, so the wake effect of the neighbouring rows is not considered. Additionally, we assume that the wind direction is parallel to the string of WTs. Optimisation problem is formulated as follows: 
\begin{align}
\mbox{max} &\quad {\sum\limits_{i=1}^n E_{k,i}},
\nonumber\\ \mbox{s.t.} &\quad P_i\leq{P_i^{max}} \quad\quad\qquad \quad\quad\quad\qquad\;\;\;\;\;\,\, \forall i\in\left\{1 \ldots n \right\}
\nonumber\\ &\quad 0\leq\beta_i^{sub}\leq\beta_i^{max} \quad\quad\quad\qquad\qquad\quad\, \forall i\in\left\{1 \ldots n \right\}
\nonumber\\ &\quad P_i^{sub}=(1-DM)P_i^{opt} \quad\quad\quad\quad\;\;\, \forall i\in\left\{1 \ldots n-1 \right\}
\nonumber\\ &\quad \omega_i^{min}\leq\omega_i^{opt}\leq\omega_i^{sub}\leq\omega_i^{max} \;\;\qquad\;\;\, \forall i\in\left\{1 \ldots n \right\}
\nonumber\\ &\quad v_{i+1}=v_{i}+k^{'}(v_1-v_{i})-kv_1C_{Ti} \; \forall i\in\left\{1 \ldots n-1 \right\}
\end{align}
\par The optimisation variables are $\omega_i$ and $\beta_i$, so the number of variables is $2n$, where $n$ is the number of WTs in the WF. In a Type III WT, $\omega^{min}$ and $\omega^{max}$ are limited by the size of the power electronic converter. Furthermore, the minimum rotor speed in the sub-optimal operation mode is further limited by the minimum rotor speed in the optimal mode. As shown in Fig.~\ref{Figure_1}, to increase the stored kinetic energy of the rotor, its speed in the sub-optimal mode should be higher than its speed in the optimal mode $\omega^{sub}\geq\omega^{opt}$. 

The resulting constrained non-linear optimisation problem is non-convex and global, having several local optima. We solved it with the pattern search algorithm using MATLAB Global Optimisation Toolbox \cite{MATLAB}. Although the optimisation problem has multiple optima, we 
were able to solve it efficiently with a good-quality solution. For a small system used in this paper, the computational efficiency was not an issue. For larger system, the optimisation can be performed off-line and solutions stored in a look-up table, which can then be readily
implemented in a control algorithm.

\section{Case study}
\par
\begin{figure}
\centering
\includegraphics [width=8.7cm, height=3.8cm]{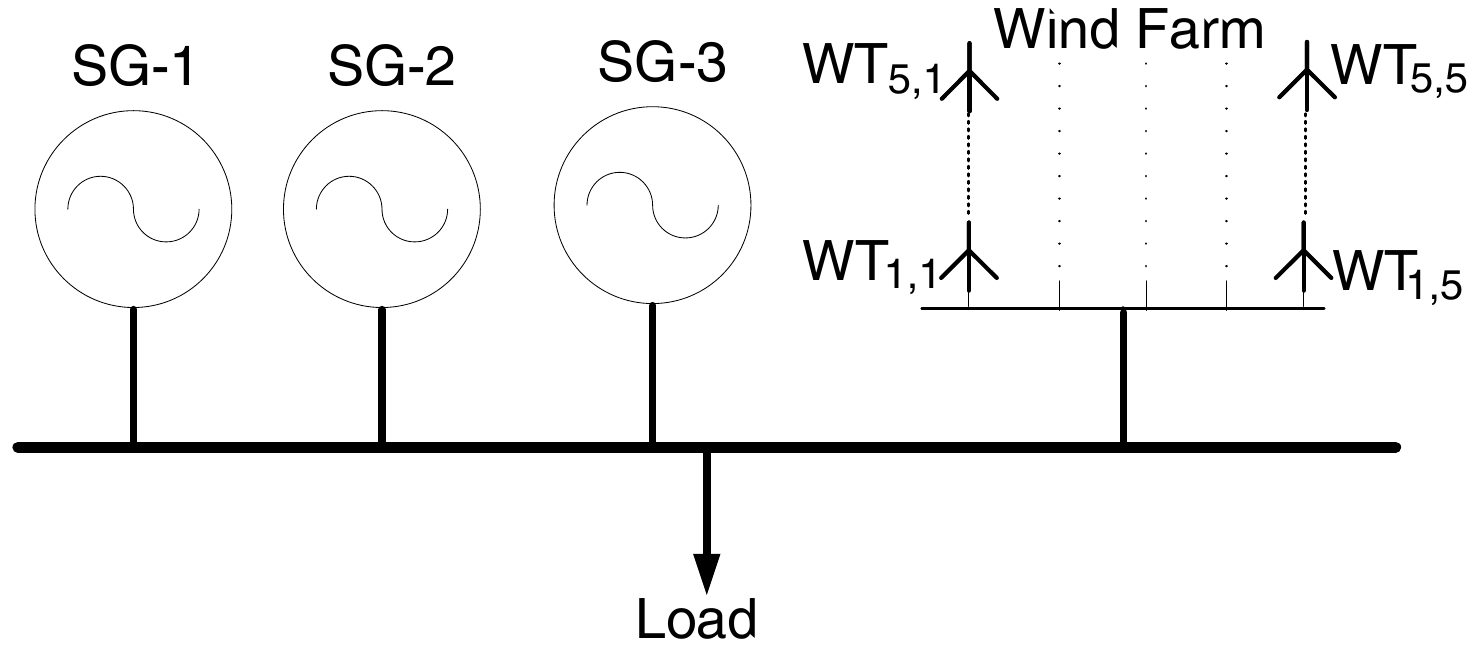} 
\caption{One-line diagram of the power system for simulation.}
\label{Figure_6}
\end{figure}

 We test the proposed strategy on a small test system with high penetration of wind generation shown in Fig.~\ref{Figure_6}. The system consists of three SGs, a WF with 5 rows of WTs, each consisting of five \SI{5}{\mega\watt} Type III WTs, and a \SI{130}{\mega\watt} constant load. We assume that the distance between the rows is large enough so we can ignore the wake interaction between the neighbouring rows, which enables us to perform the optimisation for each row independently. 

We consider a \SI{5}{\mega\watt} NREL reference WT \cite{Jonkman2009}. Using the pattern-search algorithm  with the above constraints, we can optimise the total kinetic energy of a WF with $n$ WTs using the combination of $\omega$ and $\beta$ in Zones 1 and 2. In Zones 3 and 4, the rotor speed is at its maximum limit; therefore, the only control variable is $\beta$, so the kinetic energy cannot be increased further.   

\subsection{Optimisation results}
\par
\begin{figure}
\centering
\includegraphics [width=9cm, height=6.8cm]{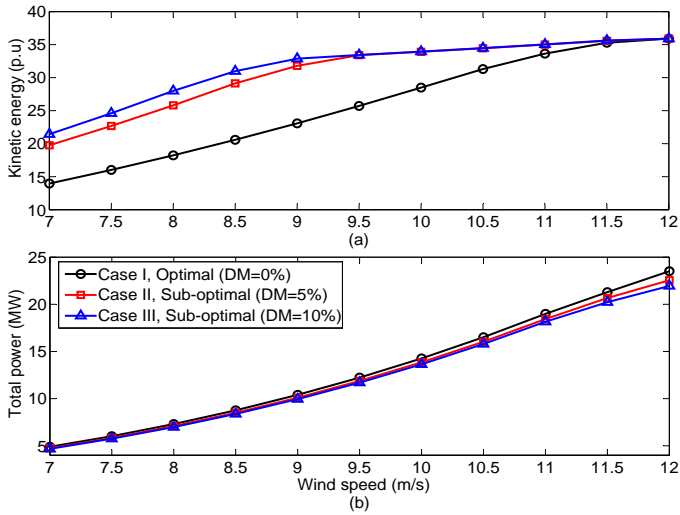} 
\caption{The optimised kinetic energy in p.u, based on the ratings of the WT, and the power of the WF with different $DMs$.}
\label{Figure_5}
\end{figure}

 The optimisation of the kinetic energy was performed for a WF with five Type III WTs in a row. $DM$ for $WT_1-WT_4$ was set to 0\%, 5\%, and 10\% for Cases I, II, and III, respectively. $WT_5$ is the last WT in the WF, so we maximise its power production without de-loading. For all three cases, the optimisation was performed for a  wind speed range of $\SI[per-mode=fraction]{7}{\meter\per\second}\leq{v}\leq\SI[per-mode=fraction]{12}{\meter\per\second}$, where the rotor speed varies between $\omega^{min}\leq\omega\leq\omega^{max}$. Fig.~\ref{Figure_4} shows the optimal rotor speeds and the pitch angles of $WT_1-WT_5$ in all cases. The optimised kinetic energy and the total power of the WF is shown in Fig.~\ref{Figure_5}. In Case I, all WTs operate in the optimal mode, and no additional kinetic energy can be stored  as shown in Figs.~\ref{Figure_4}a and 5a. In Cases II and III, the results are similar for lower wind speeds (Zones 1 and 2), and we can store considerable amount of kinetic energy by varying both the rotor speed and the pitch angle as shown in Figs.~\ref{Figure_4}b and c. Since both $C_p$ and $C_T$ are functions of the rotor speed and the pitch angle, and the operation of the down-WTs is linked to the operation of the up-WTs by the wake equation (1), we can vary both the rotor speed and the pitch angle of the WTs and not affect the total power of the WF, but still be able to store considerable amount of kinetic energy in the rotors as shown in Fig.~\ref{Figure_5}a.
\par In Cases II and III, the rotor speed and the pitch angles of $WT_1-WT_4$ follow similar trends until the rated rotor speed is reached (Zones 1 and 2). In these regions, setting the pitch angle approximately between \ang{1} and \ang{2}, and maintaining the total power of the WF close to optimal, the rotor speed increases more rapidly, thus increasing the total amount of stored kinetic energy. As shown in Fig.~\ref{Figure_5}a, we can increase the total kinetic energy of the WF by increasing the $DM$. However, once the rated rotor speed is reached, no further kinetic energy can be stored. Fig.~\ref{Figure_5}a shows that for wind speeds $v\geq\SI[per-mode=fraction]{9.5}{\meter\per\second}$ the rotor speeds are at their respective maximal limits in both cases. Therefore, no additional kinetic energy can be stored in Case III, and in this region the optimisation of the kinetic energy is independent of the $DM$. For wind speeds $v\leq\SI[per-mode=fraction]{9.5}{\meter\per\second}$ on the other hand, the available kinetic energy of the rotor is directly proportional to the $DM$ as shown in Fig.~\ref{Figure_5}a and Fig.~\ref{Figure_7}.

\par Observe that for lower wind speeds the cumulative output power of the WF is almost the same in all three cases. In higher wind speeds, on the other hand, as the $DM$ increases, the total output power decreases. In lower wind speeds, we can take advantage of the wake interaction to partly recuperate the power which we lose by de-loading. For higher wind speeds, we are unable to do this because the power that we lose  by de-loading of the up-WTs is comparatively higher, which negatively impacts the efficiency of the WF, as shown in Fig.~\ref{Figure_5}b. A better strategy for higher wind speeds would be to de-load fewer WTs with a higher DM, which would allow down-WTs not only to recuperate the power which we lose by the de-loading, but would also possibly improve the overall efficiency of the WF. Indeed, it has been shown in \cite{Madjidian2011} that in a WF with 10 WTs, de-loading the first WT by 10\% improves the row efficiency of the WF by 3\%.

\subsection{Time domain simulations}

\par To validate the optimisation results, we use time-domain simulation to simulate a frequency disturbance in the simple test system. Unlike voltage, frequency is a global variable so we can use a copper-plate model. We assume that the free wind speed $v_1$ reaching $WT_{1,1}-WT_{1,5}$ is uniform and $v_1=\SI[per-mode=fraction]{8}{\meter\per\second}$, which results in an overall output power of \SI{39.5}{\mega\watt} for the WF in all three cases. The droop is set to 1\% for all WTs. The characteristics and the settings of the SGs  are given in Table~\ref{Table_1}.

\begin{table}
\centering
\caption{Settings and Characteristics of Synchronous Generators}
\label{Table_1}
\begin{tabular}{cccccc}
\hline\hline
SG & Type  & \begin{tabular}[c]{@{}c@{}}$P_{out}$\\ (MW)\end{tabular} & \begin{tabular}[c]{@{}c@{}}$P_{max}$\\ (MW)\end{tabular} & \begin{tabular}[c]{@{}c@{}}Governor\\ Type\cite{ieeepes2013}\end{tabular} & \begin{tabular}[c]{@{}c@{}}Droop \end{tabular} \\ \hline
1  & Steam & 25.5 & 45 & IEEEG1  & 20\% \\ 
2  & Gas & 45 & 50 & TGOV1 & 5\% \\ 
3  & Gas & 20 & 25 & TGOV1& 5\%  \\ \hline\hline
\end{tabular}
\end{table}
\par
\begin{figure}
\centering
\includegraphics [width=9cm, height=3.8cm]{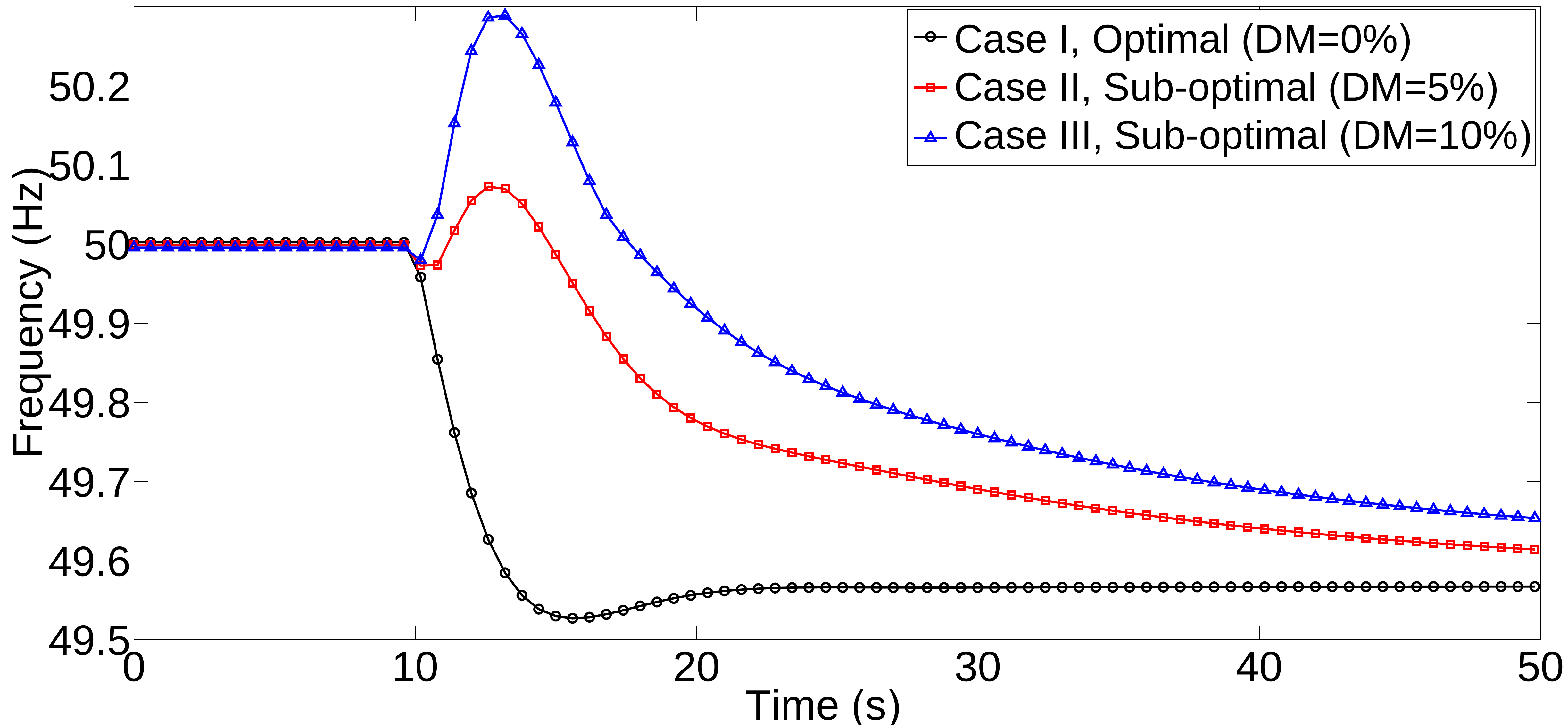} 
\caption{System frequency during generation loss.}
\label{Figure_7}
\end{figure}

\begin{figure}
\centering
\includegraphics [width=9cm, height=7.7cm]{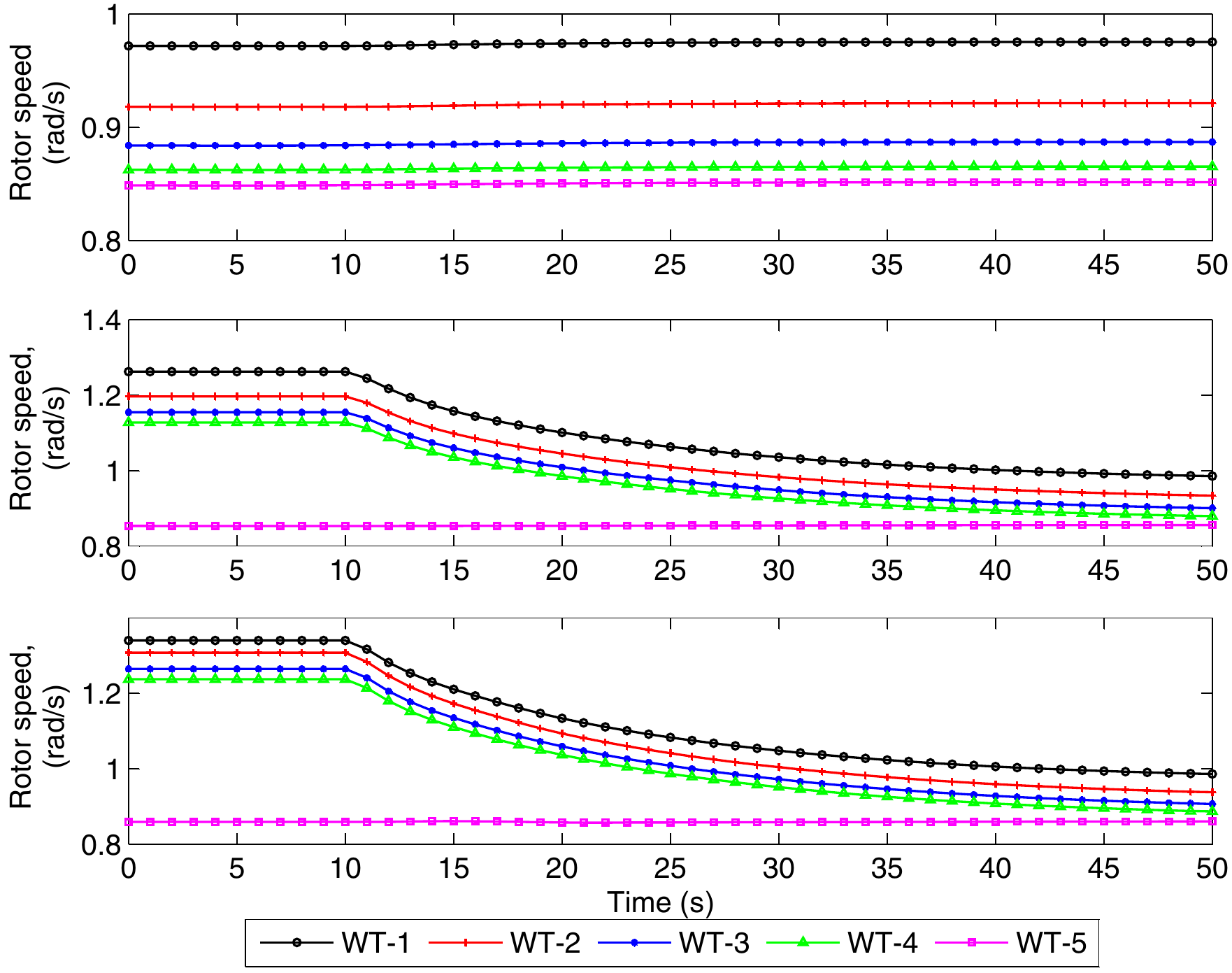} 
\caption{Type III WTs rotor speed during generation loss a) $DM=0\%$ (Case I), b) $DM=5\%$ (Case II), c) $DM=10\%$ (Case III).}
\label{Figure_8}
\end{figure}

 Simulation was performed for the three cases described in Section V.A. For all cases, we disconnect SG-3 from the system at $t=\SI{10}{\second}$. Meanwhile, in Cases II and III, we change the operation mode of $WT_1-WT_4$ form sub-optimal to optimal to release the stored kinetic energy into the system. The system frequency response is shown in Fig.~\ref{Figure_7} and the rotor speeds are shown in Fig.~\ref{Figure_8}. Since all rows of WTs in the WF are identical, only the rotor speeds of the first row is shown.     

\par In Case I, all WTs  operate in the optimal mode, so there is no additional kinetic energy. Therefore, after disconnecting SG-3, the system frequency drops rapidly and  the frequency bottom of \SI{49.52}{\hertz} is achieved after \SI{6}{\second}. It takes about \SI{13}{\second} for the frequency to plateau as shown in Fig.~\ref{Figure_7}. The rotor speeds are shown in Fig.~\ref{Figure_8}a. Because the wind speed is constant, and the WTs operate in the optimal mode, the rotor speeds remain constant. 

\par In Cases II and III, considerable amount of kinetic energy is available as shown in Fig.~\ref{Figure_5}. Therefore, after disconnecting SG-3, we inject this kinetic energy to the system. As a result, for the first few seconds, the system frequency not only doesn't drop, but it also rises, as shown in Fig.~\ref{Figure_7}. Since in Case III the $DM$ is higher than in Case II, WF contribution on the primary frequency control is even higher. Fig.~\ref{Figure_7} shows that in these two cases too much of kinetic energy is released too soon, which results in a frequency over-shoot. A better strategy would be to inject this energy into the system gradually, either by using the control strategy developed in \cite{Zertekb2012}, or by reducing the rate of change of the rotor speed. Figs.~\ref{Figure_8}b and c illustrate the rotor speeds of WTs in Cases II and III. Observe that there is a significant discrepancy between the rotor speeds of $WT_1-WT_4$ in the sub-optimal and the optimal operation modes, which results in a positive contribution of the WF to the primary frequency control. On the other hand, $WT_5$ that is not de-loaded, doesn't participate in the frequency control, and has a constant rotor speed.

\section{Conclusion}

\par In this paper, we have proposed an optimised operation strategy for de-loading operation of WFs. In this strategy, we maximise the total stored kinetic energy of the WTs' rotors, which can be released into the system during a frequency dip. We do this by optimising the rotor speeds and the pitch angles of some up-WTs. In contrast to the traditional WF operation based on optimally operating every individual WT, we consider the whole WF as a single unit by taking advantage of the wake interaction within the WF. We have shown that not only a WF can provide primary frequency control, but also it can possibly deliver more power. It has been shown that by using WTs' rotor inertia, we can delay the system frequency nadir for up to \SI{30}{\second}. We have however observed that in some scenarios too much kinetic energy could be released into the system too soon, which indicates that more intricate control strategies might be needed to achieve an optimal system performance.

\bibliographystyle{IEEEtran}
\bibliography{WindInt}

\end{document}